\documentclass[12pt,a4paper,leqno]{article}
 \usepackage{amsmath}
\usepackage{amssymb}
\usepackage{amsfonts}
 \usepackage{amsthm}
\usepackage[english]{babel}
\usepackage[applemac]{inputenc}
\usepackage{mathrsfs}

\numberwithin{equation}{section}
\renewcommand{\phi}{\varphi}
\renewcommand{\tilde}{\widetilde}

 \DeclareMathOperator{\Aut}{Aut}
 \DeclareMathOperator{\Id}{Id}

\newcommand{\ch}{\mathrm{ch}}
\newcommand{\td}{\mathrm{td}}
\renewcommand{\L}{\mathscr{L}}
\renewcommand \O {\mathscr {O}}
\def\n{\underline{n}}

\newcommand{\Hom}{\mathrm{Hom}}
\newcommand{\ra}{\rightarrow}
\newcommand{\rond}{\circ}

\def \Z {{\mathbb Z}}
\def \C {\mathbb C}
\def \R {\mathbb R}
\def \PP {\C P}
\def \QQ {\mathbb Q}

\newcommand{\codim}{\mathrm{codim}}

\begin{document}

\newtheorem{theorem}{Theorem}[section]
\newtheorem{lemma}[theorem]{Lemma}
\newtheorem{proposition}[theorem]{Proposition}
\newtheorem{Cor}[theorem]{Corollary}
\newtheorem{corollaire}[theorem]{Corollary}

\def \cal {\mathcal}

\theoremstyle{definition}
\newtheorem{definition}[theorem]{Definition}
\newtheorem{example}[theorem]{Example}
\newtheorem{exer}[section]{Exercise}
\newtheorem{conj}{Conjecture}

\theoremstyle{remark}
\newtheorem{remark}[theorem]{Remark}
\newtheorem{remarks}[theorem]{Remarks}
\newtheorem{question}{Question}

\renewcommand{\thefootnote}{\fnsymbol{footnote}}

\date{}

\title{\textbf {Estimates of Characteristic numbers of real algebraic varieties}}
\author{ \sc Yves Laszlo\footnotemark[1] \qquad Claude Viterbo\footnotemark[1]}

\maketitle \footnotetext[1]{Centre de Math\'ematiques Laurent Schwartz, UMR 7640 du CNRS, Ecole Polytechnique - 91128 Palaiseau, France. \\
email:  laszlo@math.polytechnique.fr,  viterbo@math.polytechnique.fr  }
\abstract{We give some explicit bounds for the number of cobordism classes of  real algebraic manifolds of real degree less than $d$, and for the size of the sum of $\mod 2$ Betti numbers for the real form of complex manifolds of complex degree less than $d$.}
\renewcommand{\thefootnote}{\arabic{footnote}}
\section*{Introduction}\label{intro}
In complex as well as  in real algebraic geometry, it is useful to
understand what are the simplest manifolds, and try to 
 measure the complexity of algebraic manifolds.
\textit{A priori} there are many ways to define what ``simple'' manifolds should be
and for no obvious reason should they all agree. Just to mention a
few examples, among the most frequently used (see \cite{Kollar}),
there are low degree manifolds, rational manifolds, Fano and
uniruled manifolds. These manifolds have sometimes bounded Betti
numbers, some belong to finitely many deformation classes (low
degree and Fano).

Once we decide on a given notion of complexity for complex
manifolds, we may try to understand the properties of their  real part.
Are the  real part of "simple" complex manifolds
"simpler"  ? 
 
The situation is far from being understood. Among the  oldest results in real algebraic geometry, we find the problem of estimating the Betti numbers of a real algebraic manifold.

Besides Harnack's result for the maximal  number of connected components of a real algebraic curve, dating back to 1876, and later work, for example by Severi and Commessatti for surfaces, the litterature we are aware of, may be divided in two classes. 

In the first approach,  the $n$-dimensional manifold $M^n_\R$  is the set of real points of a complex manifold, $M_\C^n$. The inequality between $\mod 2$ Betti numbers \begin{equation} \tag{1} \sum_{j=0}^n b_j(M^n_\R) \leq \sum_{j=0}^{2n} b_j(M^n_\C) \end{equation}  is then given by the so-called Smith's theory. The first occurrence of this theory can be found in \cite{Smith} and \cite{R-Smith} explicitly for the number of connected components of a real algebraic set, but the same ideas can be applied to estimate higher Betti numbers, as hinted in the paper (see particularly \cite{R-Smith} page 620, and \cite{Smith}, page 509). The explicit inequality (\thetag{1}) can first be found in \cite{Floyd}  (theorem 4.4 page 146) and the modern presentation, using spectral sequences for the  $\Z /2$-equivariant cohomology, is due to  \cite{Borel} (page 55, § 4.1). 
These results  imply that whenever there are bounds on the
{\it mod $2$}
Betti numbers of the set of complex points, there are bounds on
the set of real points.

Another approach defines $M_\R^n$ by equations $F_1=0, \ldots , F_k=0$ of degree $d_1, \ldots ,d_k$,  in $\R P^{n+k}$. 
The  Oleinik inequalities for hypersurfaces, extended to the general case by Thom and Milnor   (see \cite{Oleinik-Petrowski}, \cite{Oleinik}, \cite{Thom} and \cite{Milnor}),  are all based on applying Morse theory, usually to $F=\sum_{j=1}^k F_j^2$. Thus their estimate will depend on $\delta = \max d_k$ and yield a bound of the order of $\delta^{n+k-1}$.  

Now applying the above reults, and considering the fact that a general algebraic manifold of {\bf real} degree $d$ can be realized as a real algebraic manifold in $\R P^{2n+1}$, defined by equations of degree less than $d$, we get an estimate of the type\footnote{Alternatively, if one wants to use the {\bf complex} degree, one could use the Oleinik-Thom-Milnor inequality for the complex case, and then apply Smith theory, but the estimates would then be in $d^{4n+2}$} $$\sum_{j=0}^n b_j(M_\R^n) \leq C d^{2n+1}$$

 We will show in the second part of the paper, that  bounds on the {\bf complex} degree of the complex manifold are
sufficient to yield bounds on the Betti numbers, and  these bounds are explicit and of the order of $C d^{n+1}$. 

It is also known according to Eliashberg and the second author (cf. \cite{Kharlamov, Eliashberg}, and \cite {Kollar2} for the $3$ dimensional case), 
that complex uniruled manifolds cannot have a real form with negative sectional curvature
and this means that negatively curved manifolds are "complicated": in particular 
they can not be represented by equations of low degree.

Besides this,   little is known about the following conjectures , where the word {\it simple}
has no precise meaning\footnote{However we refer to \cite{Borel-Serre, Degtyarev-Kharlamov1, Degtyarev-Kharlamov2, Degtyarev-Kharlamov3} for some examples.}

     \begin{conj}
  Let $M(\C)$ be a smooth projective manifold. Then the set of diffeomorphism types of its real forms is finite
  \end{conj}

  \begin{conj}
   {\it Simple} complex manifolds have {\it simple} real forms.
   \end{conj}

  To prove this even in the more accessible simply connected case,
one would need bounds on the Betti numbers, which follows form the
Smith-Thom inequalities, the multiplicative structure, and the
Pontriagin classes. Our goal here is to prove that for manifolds
of low degree, the Pontriagin numbers are explicitly bounded by
the degree of a real embedding, and thus the number of cobordism classes.
  If the complex cotangent bundle
is nef, this bound depends only on the degree of a complex embedding : thus the
Pontriagin numbers are also a measure of the complexity at least
in this case.
  
We are very grateful to  Slava Kharlamov, for the many suggestions, his constant encouragement, as well as for many informations on the history of the Real Algebraic Geometry. 
We are glad to thank Arnaud Beauville and Jean-Pierre Demailly for useful discussions and advice.

 \section{Main theorem}

In this paper we shall work with complex projective manifolds
$M_\C$, with an anti-holomorphic involution $\tau$. The set of
fixed points is denoted by $M_\R$. Another point of view is to
look at a smooth projective real variety $X$ defining
$$M_\C=X(\C)\text{ and } M_\R=X(\R).$$

We'll be mainly interested in the case where $M_\R$ is further
assumed to be orientable : we'll say  simply in this case that $M$
(or $M_\R$) is an orientable real algebraic manifold.

We'll use without any further comment these two points of view
(with this notation). Of course for a given $M_\C$ there may be
many non equivalent $M_\R$ : in fact as pointed out in the introduction, 
it is not even known that
there are finitely diffeomorphism, homeomorphism or homotopy types
of $M_\R$. Vice-versa, it is not known, given $M_\R$, what are the
different complex manifolds $M_\C$ having $M_\R$ as real form. But
here even the formulation of the problem is unclear, since we may
obviously blow up $M_\C$ along subvarieties invariant by $\tau$
but not intersecting $M_\R$, without changing $M_\R$. Even in the
two dimensional case, where we may ask how many minimal $M_\C$
yield the given real from: the Klein bottle is the real form of
all the odd Hirzebruch surfaces, $ {\mathbb P}({\mathcal O}\oplus
{\mathcal O }(2k+1))$ (\cite{Silhol}).

Now given a real algebraic manifold, we shall assume it is endowed
with a {\it real} very ample line bundle $H$, where {\it real}
means that $\tau^*(H)=\overline H$. Then we denote by $\deg
(M)=c_1(H)^n$. Note that if $H$ is a (non real) very ample line
bundle over $M_C$, then $L=H\otimes \tau^*(\overline{H})$ is a real very
ample line bundle, but there is no general bound of its degree
depending only on the degree of $H$.

Our main results are as follows

  \begin{theorem}\label{boundpon}
  Let $M_\R$ be an orientable real algebraic manifold of degree $d$ and dimension $n$.
  Then the Pontriagin numbers satisfy the following inequalities

  $$ \vert p_I(M_\R) \vert \leq 2^{n^2+3n}d (d+n-2)^{n}$$

 In particular the number of  cobordism classes of manifolds representable
 as real algebraic manifolds of a given real degree can be explicitly bounded.
  \end{theorem}    

Note that the existence of such a bound can be proved a priori,
 since one can prove that real  manifolds of 
  a given (real) degree fall into finitely many diffeomorphism classes, but we believe
  making these bounds explicit is useful. 
  The only result we know of, related to the above, is due to Thom, who proves in \cite{Thom2} that if $V$ is a compact real affine complete intersection (i.e. given by $k$ polynomial equations in $\R^{n+k}$) then it is an (unoriented) boundary, and thus all Stiefel-Whitney numbers are zero. 
  
 The next step in our program should then be to show that such a bound also holds for $d$ the complex degree. This would then imply that the real part of a complex manifold of given degree falls in 
 finitely many cobordism classes, a first step towards the proof of conjecture $1$. For the moment, there are only few classes for which this can be proved. Note that, as pointed out by Kharlamov,  using Torelli's theorem for K3 surfaces, it is easy to construct surfaces with compllex degree bounded and  the real degree can be arbitrarily large.

We also give some estimate on the size of the homology for the real part of a complex manifold. Note that contrary to the previous theorem, the degree used here is the complex degree, that is the minimal degree of any algebraic embedding, not necessarily real.

Let $b(X)$ represent the sum of the {\it mod} $2$ Betti numbers. 

\begin{theorem} Let $X$ be a real subvariety of the projective
space. Assume that $X_\C$ is smooth, connected of complex degree $d$ and
dimension $n$. Then, one has $b(X(\R))\leq b(X(\C))\leq 2^{{n^2}+2}d^{n+1}$.
\end{theorem}

The first inequality is just the Smith-Thom inequality. As for the second one, according to Milnor (\cite{Milnor}), a complex submanifold in ${\mathbb C} P^m$ given by equations of degree less than $d$ has $b(X_\C) \leq C_n\delta^{2m+2}$. Since if $X(\C)$ has degree $d$, it is given by polynomials equations of degree at most $d$ in $\C P^{2n+1}$ (see \cite{Mumford}),  Milnor's inequality yields $b(X_\C) \leq C_n d^{4n+4}$. 
In the real case, for a submanifold in $\R P^m$, Milnor gives $b(X_ \R) \leq C_n \delta^{m+1}$ and the same argument yields $b(X_\R) \leq C_nd^{2n+2}$. However we cannot say that we improve on Milnor's inequality, since there seems to be no clear lower bound for the maximal degree of the equations defining $X$ in $\C P^{2n+1}$ as a function of the degree of $X$. 

\section{Some lemmata on real and complex cycles} 
Since we shall deal with non-orientable situations, we shall use homology and cohomology with local systems of coefficents\footnote{We shall only use the case of  local coefficients with fibre  $ {\mathbb Z} $, and this implies that we always have $ \mathscr L \otimes \mathscr L = {\mathbb Z} $. Of course we can also work with $ {\mathbb Z} /2 {\mathbb Z} $ coefficients, but then $ \mathscr L $ is trivial and much of what follows is much simpler.}, $ \mathscr L $. Unless otherwise stated, all our manifolds are closed and connected. We denote by $H^*(X, \mathscr L )$ and $ H_*(X, \mathscr L )$ such groups, by $H_c^*(X, \mathscr L )$ its compact supported version. We denote by $o_X$ the orientation bundle of $X$. 

Moreover we have cup and cap  product map

 \begin{gather*} 
\cup : H^k(M, \mathscr L ) \otimes H^l (M, \mathscr L ') \longrightarrow H^{k+l}(M, \mathscr L \otimes \mathscr L ') \\
\cap :  H^k(M, \mathscr L ) \otimes H_l (M, \mathscr L ') \longrightarrow H_{l-k}(M, \mathscr L \otimes \mathscr L ') 
 \end{gather*} 
 
The first map is a non-degenerate pairing  for $k+l=n, \mathscr L \otimes \mathscr L'=o_M$.

 Also a map $f:V \to M$ induces maps 
 
 $$f_*: H_*(V, f^*( \mathscr L )) \longrightarrow H_*(M, \mathscr L )$$
 $$f^*: H^*(V, \mathscr L )  \longrightarrow H^*(M, f^*( \mathscr L ))$$
 
 Moreover we have a Poincaré duality map:
 
 $$ PD: H_k(M, \mathscr L ) \longrightarrow H^{n-k}(M, \mathscr L \otimes o_M)
 $$
 This implies in particular that $H_n(M, o_M)=H^0(M, {\mathbb Z} )$ is canonically isomorphic to $ {\mathbb Z} $. 
 
 Finally to a $k$-dimensional vector bundle, $\pi:E \to X$,  we associate the relative orientation bundle $o_{E/X}$, as the set of orientations of the fiber, namely $o_E\otimes \pi^*(o_X^{-1})$. The Thom isomorphism theorem is  a map $$T:H^*(X, \mathscr L ) \longrightarrow H_c^{*+k}(E, \mathscr L \otimes o_{E/X})$$ and this is given by $T(u)= \pi^*(u)\cup T_E$ where  $T_E\in H_{c}^k(E, o_{E/X})$. 
 The class $T_E$ is uniquely defined as the preimage of $1$ by the canonical isomorphism $ H^{k}_c(E, o_{E/X}) \to {\mathbb Z} $. Naturality implies that for every $x$ in $X$,  $T_E$  is uniquely determined by the property that it restricts on $H^k_c(E_x, o_{E_x})$ to the canonical generator.

 From this it easily follows\footnote{Note that if $\pi, \pi'$ are the projections of $E\oplus E'$ on $E,E'$, then $\pi^*(T_E)$ has compact support only in the $E$ direction, and $(\pi')^*(T_{E'})$ in the $E'$ direction, but then $\pi^*(T_E) \cup (\pi')^*(T_{E'})$ has compact support, and we denote this class by $T_E\cup T_{E'}$. }, using the canonical isomorphism $\pi^*o_{E/X}\otimes\pi'*(o_{E'/X}\simeq o_{E\oplus E'/X}$,
 that $T_{E\oplus E'}=T_E\cup T_{E'}$, since we just have to check this over a point. 

Finally note that Poincaré duality on $E$ is an isomorphism 
$$H_{n-d}(E, \mathscr L ) \to H^{k+d}_c(E, \mathscr L \otimes o_M)$$
but since $M$ is a retract of $E$, and $$H_{n-d}(E, \mathscr L ) \simeq H_{n-d}(M, \mathscr L ) \simeq H^d(M, \mathscr L \otimes o_M)$$ the second isomorphism being again Poincaré duality, 
we get an isomorphism $$H^d(M, \mathscr L \otimes o_M) \to  H_c^{k+d}(E, \mathscr L \otimes o_M)$$
and so the Thom isomorphism is induced by Poincaré duality. 

 Now if $ \mathscr L $ is a local system and we have
 a map $j:V\to M$ and $ j^*( \mathscr L )\otimes o_V\simeq {\mathbb Z} $ (or equivalently  $ j^*( \mathscr L ) \simeq o_V$), then we say that $V$ is $ \mathscr L $-coorientable. If moreover one of these two isomorphism is given, we say that $V$ is $\mathscr L$-cooriented. 
 In the case $ \mathscr L = o_M$ we just say that $V$ is coorientable (resp. cooriented).  
 
 Now if $V$ is  $k$-dimensional, connected and $ \mathscr L $-cooriented, we may associate to a $V$ a class $[V]\in H_k(M, \mathscr L)$ as follows. By Poincaré duality we have an isomorphism $$H_k(V, j^*( \mathscr L )) \to H^0(V, j^*( \mathscr L ) \otimes o_V) \simeq {\mathbb Z} $$ and the last isomorphism is canonically given  by the $\mathscr L$- coorientation of $V$. In the non connected case, we just do the same for each connceted component. 
 
 Then we have $j_*: H_k(V, j^*( \mathscr L )) \to H_k(M, \mathscr L )$ and we still denote by $[V]$ the class $j_*([V])\in H_k(M, \mathscr L )$. 
 
The submanifold $V$ also has a fundamental class $\mu_V$ in $H^{n-k}(M, \mathscr L\otimes o_M)$  defined as the extension of the Thom class of the normal bundle of $V$ in $M$. It is the image of $[V]\in H_k(M, \mathscr L )$ by Poincaré duality. 


Note that if $\beta \in H^k(M,  \mathscr L ^*)$ then $<\beta, [V]> \in H_0(M, {\mathbb Z} )\simeq {\mathbb Z} $ is given by $\mu_V\cup \beta \in H^n(M,  \mathscr L \otimes \mathscr L ^* \otimes o_M)\simeq H^n(M, o_M) \simeq {\mathbb Z} $.

 \begin{lemma} 
 Let $V,W$ be $ \mathscr L $-cooriented and $ \mathscr L '$-cooriented submanifolds of $M$. Assume $V,W$ are in general position, and  set $Z=V\cap W$. Then $Z$ is $\L\otimes\L'\o_M$ coriented and if $\mu_V \in H^{k}(M, \mathscr L \otimes o_M)$ and  $\mu_W \in H^{l}(M, \mathscr L '\otimes o_M)$ then $\mu_V\cup \mu_W = \mu_Z \in H^{k+l}(M, \mathscr L \otimes \mathscr L ')$
 \end{lemma}

 \begin{proof} 
 
 The class $\mu_V$ is related to the Thom class as follows: if $V$ is an $ \mathscr L $ cooriented submanifold of $M$, $N(V,M)$ its normal bundle, and $T_{(V,M)} \in H^{n-k}_c(N(V,M),  \otimes o_{N(V,M)})$ its Thom class, then the extension of $T_V$ to $H^{n-k}(M,  \otimes o_{M})$ is precisely $\mu_V$. 
 
  Moreover by naturality of the Thom class, if $i:Z \to V$,  $E$ is a bundle over $V$, and we denote by $\tilde i$ the natural map $i^*(E)\to E$, that $T_{i^*(E)}=({\tilde i})^*(T_E)$.   
  
 Now let $j_V:V \cap W \to V, j_W: V\cap W \to W$ be the inclusions, since $j_V^*(N(V,M))=N(V\cap W,W), 
 j_V^*(N(W,M))=N(V\cap W,V)$, and of course by general position $$N(V\cap W, M)=N(V\cap W, W) \oplus N(V\cap W,V)$$ that is 
  $$N(V\cap W, M)=j_W^*N(W, M) \oplus j_V^*(N(V,M))$$

 and this induces an isomorphism $o_{N(V\cap W,M)/V\cap W}\simeq j_V*o_{N(V,M)/V\cap W}\otimes j_W*o_{N(W,M)/V\cap W}$ and we get we  get that $T_{N(V\cap W,M)}=T_{N(W,M)}\cup T_{N(V,M)}$,   that is

 $$T_{V\cap W}=T_V\cup T_W$$

 Thus   if $V,W$ are $ \mathscr L $ and $ \mathscr L '$-coorientable submanifolds, in general position, we have
  
  $$T_{V\cap W}=T_V \cup T_W$$ 

this implies
  
 $$ \mu_{V\cap W} = \mu _V \cup \mu_W$$
 
 \end{proof} 
 
 In particular  if $V$ is oriented (and therfore $\Z$-cooriented) and $W$ cooriented, we have that $\mu_V \in H^*(M, o_M)$, $\mu_W\in H^*(M,\Z)$ then $Z$ is oriented and dual to $\mu_V\cup \mu_W \in H^*(M,o_M)$.

\begin{remark} \label{singsing}

This extends to the case where $V,W$ are singular submanifolds with codimension $2$ singular locus. 
Indeed in this case the fundamental class is also well defined, and  provided $V,W$ are in general position, $V\cap W$ is also a singular manifold with codimension $2$ singular locus and the above result also applies. (see \cite{Borel-Haefliger} for more details). 
\end{remark} 
 
Let's now consider a smooth real projective
manifold $M$ of dimension $n$.

Let $V,W$ be two real submanifolds of $M$, so that we have
$V_{\C}, W_\C\subset M_\C$ and $V_\R,W_\R \subset M_\R$. Note that $V_\C, W_\C, M_\C$ are always orientable

Let us denote by $`` \cdot "$ the algebraic intersection of cycles, and
assume $V_\C, W_\C$ are in general position and have complementary dimension. For example if $V_\C, W_\C$ are
zero sets of sections of very ample bundles, we can always move
$V,W$ so that this assumption is satisfied. 

 \begin{lemma}\label{inter-bound}
Assume $V_\R$ is $ \mathscr L $-coorientable, $W_\R$ is $ \mathscr L^* $-orientable, and $V_\R , W_\R $ are in general position. Then  $$ \vert V_\R \cdot W_\R \vert \leq V_{\mathbb C} \cdot W_\C$$
 \end{lemma}
 \begin{proof}The assumption is only needed so that the number $V_\R \cdot W_\R$ is well defined. 
Indeed assume the submanifolds are in general position. Then
$\vert V_\R \cdot W_\R \vert $ counts with sign the number of
intersection points of $V_\R$ and $W_\R$. But these are contained
in the set of intersection points of $V_\C$ and $W_\C$ while
$V_{\mathbb C} \cdot W_\C $ counts all intersection points with
positive sign. The inequality is now obvious. \end{proof}

  \begin{remark} 
  
  The same proof shows that if $V^j_ {\mathbb R} $ are $ \mathscr L _j$-coorientable, with $$ \mathscr L _1\otimes \mathscr L _2 \otimes  \ldots \otimes \mathscr L _k \otimes (o_M)^k= {\mathbb Z} $$ then  have
  
  $$ \vert V^1_ {\mathbb R} \cdot  V^2_ {\mathbb R} \ldots V^k_ {\mathbb R} \vert  \leq 
  V^1_ {\mathbb C} \cdot  V^2_ {\mathbb C} \ldots V^k_ {\mathbb C}$$
  \end{remark}  
  \bigskip
Let $f$ be a "real" holomorphic map from  a real manifold $V$ to
$M$, that is $f:V_\C \to M_\C$ commutes with the anti-holomorphic involution, hence
sends $V_\R$ to $M_\R$. Then $f$ sends $[V_\C]$ to a class
$(f_\C)_*([V_\C])$, and $[V_\R]$ to $(f_\R)_*([V_\R])$. We claim
that the above result also implies that

 \begin{proposition}  \label{f-maps}Assume $V_\R$ is $f_\R^* (\mathscr L )$-coorientable, $W_\R$ is $f_\R^* (\mathscr L^* )$-orientable.  If $f(V_\C) \cap g(W_\C)$ is finite, we have
$$ \vert (f_\R)_*([V_\R])\cdot (g_\R)_*([W_\R]) \vert \leq (f_\C)_*([V_\C])\cdot (g_\C)_*([W_\C])$$
 \end{proposition}
 \begin{proof}
With the same argument as above, we just have to prove that the multiplicity of a complex intersection is always larger than the multiplicity of the real intersection. But this is a local statement, and we can then of course move locally $V$ and $W$ to put them in general position, in which case the statement is obvious.  \end{proof}
 \begin{remark}
If $V\cap W$ contain a component of positive dimension, the above inequality may fail, since $V_\C \cdot W_\C$ can be negative.
 One could however hope for an inequality with $V_{\mathbb C} \cdot W_\C$
 replaced by $ \vert V_{\mathbb C} \cdot W_\C \vert $.
 Note that the main example, is given by a  "real" holomorphic bundle $E$
  over $V_\C$ and $W=V$.  Then $V\cdot W$ denoted by $\chi(E)$ is the number of zeros of a generic smooth section of $E$ counted with sign. In the complex case this coincides with $c_n(E)$ ($n=\dim (E)$). In both real and complex case, if $E=TV$ it co\"\i ncides with $$\sum_{j=0}^{\dim (V)}(-1)^j\dim H^j(V, {\mathbb Z} ).$$
 Then $V_\C\cap W_\C=\chi (E_\C)$, and  $V_\R\cap W_\R=\chi (E_\R)$. 
  Thus we could expect $$\vert \chi (E_\R) \vert \leq \vert \chi (E_\C) \vert $$ 
However there is the following counterexample,
  suggested by Slava Kharlamov: \newline let $M_\C$ be a four dimensional complex manifold such that $M_\R$ is orientable, and  $ \vert \chi(M_\R) \vert $ is arbitrarily large. 

Let $P$ be the blow up of $M$ along the curve $C$. The formula $\chi (A\cup B)= \chi (A) + \chi (B) -\chi (A\cap B)$ follows from the Mayer-Vietoris exact sequence and implies that  $$\chi (P_\C)= \chi (M_\C)+(\nu-1)\chi (C_\C)$$ where $\nu = \text{codim} (C)=3$, and $$\chi (P_\R) = \chi (M_\R) + 2\chi (C_\R) $$

 Since $C_\R$ is a union of circles, $\chi (C_\R)=0$ and thus  $\chi (P_\R)=\chi(M_\R)$, while by a suitable choice of  $C$, we may arrange that $\chi (P_\C)$ is small. Indeed, take $M_\C$ to be 
the product of a surface with real part of arbitrarily large Euler characteristic, and the surface $$S=\{(x,y,z,t)\in \C P^3 \mid x^2+y^2+z^2+t^2=0\}$$
 Then $S_\R$ is a sphere, $\chi (S_\R)=2$ and $S_\C$ contains a curve of genus $4$ by taking the intersection of $S_\C$ with a generic hypersurface of degree $3$ yields a curve of genus $4$, thus $\chi(C)=-6$. Since we may move $C$, after blowing-up we may find a new curve of genus $4$ in the blown-up manifold, and thus we may repeat the blow up, eventually get a manifold $P$ such that $ \vert \chi (P_\C) \vert $ is between $0$ and $5$, while $ \vert \chi(P_\R) \vert  = \vert \chi (M_\R) \vert $ is large. 

Thus we see that $ \vert \chi (P_\R) \vert = \vert \chi (TP_\R) \vert  $ can be arbitrarily large,  depending on the choice of the first surface, while $\vert \chi (P_\C) \vert = \vert \chi (TP_\C) \vert  $ is between $0$ and $5$, and $TP_\R$ and $TP_\C$ are the counter-examples we were looking for.
 \end{remark}

 \begin{lemma}\label{intersec}
 Let $X, Y$  be respectively  $ \mathscr L $-cooriented and  $ \mathscr L '$-cooriented submanifolds of $M$, and $Z$ be the clean intersection of $X$ and $Y$, that is $T_ZX\cap T_ZY=TZ$. Then if  $ \mu _X, \mu _Y $ are the cohomology classes associated to $X,Y$ in $H^*(M, \mathscr L \otimes o_M)$ and $H^*(M, \mathscr L '\otimes o_M )$ we have $ \mu _X \cup \mu _Y= \mu _N \in H^*(M, \mathscr L \otimes \mathscr L ') $ where $N$ is the zero set of a generic section of $\nu Z$, where $\nu Z=T_ZM/(T_ZX+T_ZY)$
 \end{lemma}

  \begin{proof}
  Let us take local coordinates in a neighborhood of $Z$, associated to the decomposition $T_ZM=TZ\oplus \nu_X\oplus \nu_Y\oplus \nu Z$, where $\nu_X$ (resp. $\nu_Y$)
  is the normal bundle of $Z$ in $X$ (resp $Y$). Then $X$ is parameterized by the total space of $\nu_X$: $(z,\xi_X)\to (z,\xi_X,0,0)$,  $Y$ by the total space of $\nu_Y$: $(z,\xi_Y)\to (z,0,\xi_Y,0)$ .  A generic perturbation $\tilde X$ of $X$ will be given by
  $(z,\xi_X)\to (z,\xi_X,0,  \varepsilon (z,\xi_X))$, where $ \varepsilon (z,\xi_X)$ vanishes for $ \vert \xi_X \vert \geq \delta$. Then the intersection $\tilde X \cap Y$ is given by the equation $ \varepsilon (z,0)=0$, where $ \varepsilon $ is a section of $\nu Z$. This intersection is transverse provided $\frac{\partial \varepsilon }{\partial z}(z,0)$ has maximal rank at points where $ \varepsilon (z,0)=0$, and $\tilde X \cap Y$ is equal to the zero set of the section $ \varepsilon $.  We let the reader check that the coorientations match. 
   \end{proof}

In the sequel, given a submanifold $V$ of $M$, we denote by $V\cdot V$ the intersection of $V$ with the image of $V$ by a small generic isotopy. If $V$ is $ \mathscr L $-coorientable, then $V\cdot V$ is $ \mathscr L ^{\otimes 2}\otimes o_M$ -coorientable. The class $\mu_V$ is in $H^*(M, \mathscr L \otimes o_M)$ , and $V\cdot V$ is Poincaré dual to $\mu_V^2 \in H^*(M, \mathscr L ^{\otimes 2})$. 

 We now denote by $c_V,c_W$  the codimension of $V_\R,W_\R$ in $M_\R$, so that these are also the complex codimensions of $V_\C,W_\C$ in $M_\C$. Let  $ \gamma_V \in H^{2c_V}(M_\C, {\mathbb Z} )$ be the  class Poincar{\'e} dual to $V_\C$. Assume $V$ is coorientable and  $ \rho _V \in H^{c_V}(M_\R, o_V)$
be the  class Poincar{\'e} dual to $V_\R$. Then $[V_\R]\cdot [V_\R] \in $ defined as $(T_{\nu
V_\R}\cap [V_\R]$) is dual to  $\rho_V^2\in H^{2c_V}(M_\R, {\mathbb Z} )$. 
 
Let $i$ be the inclusion of $M_\R$ into $M_\C$,  we want to compute $i^*(\gamma_V)$ in terms of
$\rho_V$.

 \begin{proposition}\label{basic-square}
 We have $$i^*(\gamma_V)=\rho_V^2$$ \end{proposition}

 \begin{proof}

Applying lemma \ref{intersec} to $X=V_\C$ and $Y=M_\R$, with
$Z=V_\R$, we have that if  $\nu (V_\R)$ is the normal bundle of
$V_\R$ in $M_\R$, it is easy to see that $J\nu(V_\R)$ is normal to
the space $T(V_\C)+ T(M_\R)$ in $TM_\C$ , and thus, a perturbation
of $V_\R$ will intersect $V_\C$ along the zeros of a section of
the bundle $J\nu(V_\R)$ over $V_\R$. Since this bundle is
isomorphic to $\nu(V_\R)$, this is the same as the zero set of a
section of the normal bundle, and this coincides with $V_\R\cdot
V_\R$.
 \end{proof}

 \begin{proposition}
If $V_\R$ is $ \mathscr L $-coorientable, we have $$(f_\C)_*([V_\C])\cap M_\R = (f_\R)_*([V_\R])\cdot (f_\R)_*([V_\R])$$

  \end{proposition}

  \begin{remark}
  Note that this is still true for $V_\C$ a variety with singularities of
  codimension at least $2$, so that the same holds for $V_\R$.
  This does not hold if the singularities have codimension $1$:
  in this case we cannot even guarantee that $V_\R$ is a cycle.
  \end{remark}
\begin{remark}
   Given the map $f:X\to Y$ the normal bundle is the quotient $f^*(TY)/TX$
   \end{remark}
    \begin{proof}
   If $f$ is an embedding, this follows from the previous proposition. Consider now the case where $V\subset E$ and $f=\pi_{\mid V}$ where $\pi:E\to B$ is a real holomorphic projection (so that $\pi^{-1}(B_\R)=E_\R$). Let $Z\subset B_\R$  be a cocycle and consider

  \begin{gather*} (\pi_\C)_*([V_\C]) \cdot B_\R \cdot Z= V_\C \cdot  (\pi_\C)^{-1}(B_\R \cap Z)=(V_\C \cdot E_\R) \cap \pi^{-1}(Z) =\\ (V_\R\cdot V_\R) \cap  \pi^{-1}(Z) = \pi_*(V_\R)\cdot \pi_*(V_\R) \cap Z
    \end{gather*}
    Since this holds for any cycle, we have that

    $$ (\pi_\C)_*([V_\C]) \cdot B_\R = \pi_*(V_\R)\cdot \pi_*(V_\R) $$

    Now the general case follows from the fact that a general real holomorphic map $f$ is the composition of an embedding $\tilde f : V \to V \times M$ given by $x \to (x,f(x))$ and the restriction of the projection $V\times M \to M$.
     \end{proof}

     \begin{remark}
 After proving the above result, we realized that it can be traced back to \cite{AK93},
 (cf. Theorem A, (b), page 311) where this was proved mod $2$, and is stated there as
 the identification
     $${\bar H}^{2k}_{{\mathbb C}-alg}(M_\R)=\{ \alpha ^2 \mid \alpha \in H^k_{alg}(M_\R)\}$$
  where ${\bar H}^{2k}_{{\mathbb C}-alg}(M_\R)$ is the set of pull-backs
  in $M_\R$ of the Poincar{\'e} dual classes in $M_\C$ of real algebraic submanifolds
  \footnote{This should not be confused with   ${H}^{2k}_{{\mathbb C}-alg}(M_\R)$
  obtained by pulling back the classes dual to any complex submanifold of $M_\C$.}, and $H^k_{alg}(M_\R)$ is the set of classes Poincaré dual to a real algebraic set.
      The proof in \cite{AK93} is more algebraic, but proves also that any
      "square of an algebraic class" is induced by a complex class mod $2$.
      However we  really need the integral coefficient case in order
      to get our estimates on the Pontriagin classes (mod $2$ estimates would be useless here).

One should always be careful that  $[V_\R]\cdot [V_\R]$ is only a square if $V_\R$ is orientable, as we see from the following example. 

     \end{remark}
  \begin{example}
1)   Let us consider the inclusion $\R P^n \to {\mathbb C} P^n$,
and the pull-back of the generator $u \in {H^2}(\C P^n)$ to
${H^2}(\R P^n)$. Remember that $H^q(\R P^n)=\Z/2\Z$ for $q$ even
different from $0,n$ and is equal to $\Z$ for $q=0$ and for $q=n$
if $n$ is even. The pull-back of $u$ is then equal to the
generator of ${H^2}(\R P^n)$. This generator can be identified with
the class $a^2 \in {H^2}(\R P^n, \Z/2 \Z)$, where $a$ is the
generator of $H^1(\R P^n, \Z/2\Z)$. Even though this class is a
square in ${H^2}(\R P^n, \Z/2 \Z)$, it isn't a square in ${H^2}(\R
P^n, \Z)$

2) Let $Q$ be the manifold of degree $2n-2$, $$Q_\C=\{(z_0, ....,
z_{2n-1}) \mid z_0^2-z_{2j}^2-z_{2j+1} ^2=0, j=1..n-1\}$$
  in ${\mathbb C} P^{2n}$.  Then $Q_\R=T^{n-1}$. The pull back of the
  hyperplane class, dual to
  $\C P^{2n-1}$ is given by the square of $T^{n-1}\cap {\mathbb R} P^{2n-1}=T^{n-2} \cup T^ {n-2}  $
  \end{example}

 For $I=(a_1,...,a_q)$, let $$S_I(k) = \{L \mid \dim (L\cap k^{n-k+i+a_i})\geq i\} \subset G_q(k^n)$$
 This is the Schubert cycle associated to the multi-index $I$. If $ \vert I \vert =\sum_{l=1}^q a_j$
the cycle has   dimension $ \vert I \vert $.
 If $1_r$ denotes a sequence of $r$ ones, $S_{1_r}(\C)$ is Poincar{\'e} dual to the $r$-th Chern class.
More generally it is important to notice that the $S_I$ are submanifolds with codimension $2$ singularities (cf.\cite{Pontryagin}), so that according to remark \ref{singsing} the previous results apply. 
   \begin{remark}
   The cycle ${\tilde S}_{1_r}(\R)$, lift of $S_{1_r}(\R)$ to the Grassmannian of oriented $k$-subspaces, ${\tilde G}_k(\R^n)$, is orientable in ${\tilde G}_k(\R^n)$ provided $r,k$ are both even (\cite{Pontryagin}).
Note that ${\tilde G}_k(\R^n)$ being simply connected is
orientable, and orientability and coorientability coincide in this
case.

Note also that   the ${\tilde S}_I(\R)$ are invariant by the
involution $\tau$ of  ${\tilde G}_k(\R^n)$ that sends a space to
the same one with opposite orientation.  If this involution is
orientation preserving, then ${G}_k(\R^n)$ will be orientable, and
this happens precisely for $k$ even.   However, the involution
changes the orientation of $S_I(\R)$, except when $I$ corresponds
to the Euler class. So in the sequel, squares are never real
squares, and in particular, one should not conclude that our
Pontriagin classes are squares (see example 3.13) !

\end{remark}

  Now let $S$ be a positive linear combination of Schubert cycles, and $c_S$ be the Poincaré dual class
  of $S( \C )$, $o_S$ be either the trivial bundle or the unique non-trivial local coefficients . We denote by $j$ the inclusion $j: G_k (  {\mathbb R}) \to G_k ( {\mathbb C} )$. Then $p_S=j^*(c_S)$.  Note that if $c_S = \sum_{\alpha \in \{1,n\}^r } a_\alpha c_1^{ \alpha_1} \cdots c_r ^{ \alpha _r}$ (the $a_\alpha$ need not be all non-negative) then \footnote{remember that only the $p_{4k}$ are non-zero} $p_S= \sum_{\alpha \in \{1,n\}^r } a_\alpha p_2^{ \alpha_1} \cdots p_{2r} ^{ \alpha _r}$. We denote by $o$ the orientation bundle of $G_k( {\mathbb R} ^n)$. 
  
 \begin{corollaire}
 Let $S(\C)$ be a positive linear combination of Schubert cycles in $G_q(\C^n)$, and $S(\R)$ the analogous cycle in $G_q(\R^n)$.
 Consider the inclusion $i:G_q(\R^n)\to G_q(\C^n)$ and $\gamma_S \in H^{ 2 \codim(S) }(G_q(\C^n), {\mathbb Z} )$
 and $\rho_S\in H^{ \codim(S) }(G_q(\R^n), o)$
 be the cohomology classes Poincar{\'e} dual to $S( {\mathbb C} )$ and $S( {\mathbb R} )$. Then $$j^*(\gamma_S)=\rho_S^2$$

In other words, the cohomology class  Poincaré dual to $S(\C)$ in $G_k( \C^n)$ induces on $G_k( {\mathbb R} ^n)$ a cohomology class  Poincaré dual of $S(\R)\cdot S(\R)$.
 \end{corollaire}
 
   \begin{proof} 
   Follows immediately from prop \ref{basic-square}.  Note that we assume here that $S( {\mathbb R} )$ is coorientable, and according to Pontriagin, this is the case exactly when $S( {\mathbb R} )$ represents a non-zero Pontriagin class. Had we taken $ {\mathbb Z} /2 {\mathbb Z} $ coefficients, we could also take the pull-backs of the odd Chern classes $c_{2k+1}$. 
   
   \end{proof}

    Let now $E\to M$ be a "real" holomorphic bundle, that is, there is  a real
     bundle $E_\R$
      over $M_\R$  such that $$E_{\mid M_\R}=(E_\R\otimes \C)$$ and the anti-holomorphic
      involution
      $\tau$ lifts to $E_\C$ in such a way $E_\R$ is the set of fixed
      points. In a more algebraic setting, $E$ comes from an
      algebraic locally free sheaf on real algebraic variety $X$.
     Assume now that  $E_\C\to M_\C$ is generated by sections,
     that is there are sections $s_1,..., s_m$ such that for each
     point $x\in M_\C$, $s_1(z), ... , s_m(z)$ generate $E_z$.
     A holomorphic section is said to be real, if $s(\tau (z))= \bar s (z)$.
     Clearly, if $E$ is generated by sections, it is generated by real
     holomorphic sections. Indeed, the sections
     $\Re (s_j)=\frac{s_j(z)+\bar s_j(\tau (z))}{2}, \Im (s_1)=\frac{s_j(z)-\bar s_j(\tau (z))} {2i}$ are real holomorphic and generate $E_\C$.

     Now the real holomorphic sections induce a real holomorphic map  $\tau_\C: M_\C \to G_k(\C^m)$, such that its restriction $\tau_\R $ to $M_\R$ has its image in $G_k(\R^m)$. Moreover, $(\tau_\C)^*(U_\C)=E_\C, (\tau_\R)^*(U_\R)=E_\R$. 

 This last remark together with  lemma \ref{inter-bound} implies
  \begin{theorem}
 Let $E_\C$ be a bundle generated by its sections over $M_\C$ and assume $E_\C$ to be real.
 Let  $J=(j_1,...,j_q)$ be such that $4\sum_{t=1}^q j_t = \dim M$ then we have the  inequality

 $$ \left\vert \left <\prod_{t=1}^q p_{4j_t}(E_\R), [M_\R]\right > \right\vert
 \leq \left <\prod_{t=1}^qc_{2j_t}^2(E_\C),[M_\C]\right > $$

     that we write for short $$\left\vert \left <p_{J}(E_\R), [M_\R]\right > \right\vert \leq \left
     <c_{J}(E_\C)^2,[M_\C]\right > $$
 \end{theorem}
 \begin{proof}
    
     With $i$ denoting the inclusion $i: M_ {\mathbb R} \to M_ {\mathbb C} $ we have 
     $\tau_ {\mathbb C} \circ i = j \circ \tau_ {\mathbb R}$. Then if $S$ is a cycle  as in the previous 
     corollary we have 
     \begin{gather*} \left < p_S(E_\R ), M_\R \right > = \left < \tau_\R ^* j^*(PD(S( {\mathbb C} )), M_\R \right > =\\ \left < \tau_\R^*(PD(S( {\mathbb R} ) \cdot PD( S( {\mathbb R} ))), M_\R\right > =\left <  \tau_\R^*(PD(S( {\mathbb R} )) \cdot \tau_\R^*(PD(S( {\mathbb R} )) , M_\R \right > \end{gather*} 
     
     According to proposition \ref{f-maps}, we have the inequality
     \begin{gather*}  \vert  \left < \tau_\R^*(PD(S( {\mathbb R} )) \cdot \tau_\R^*(PD(S( {\mathbb R} )), M_\R \right >  \vert \leq \\ 
    \left <  \tau_\C^*(PD(S( {\mathbb C} )) \cdot \tau_\C^*(PD(S( {\mathbb C} )), M_\C \right >  =
    \left < c_S^2, M_\C \right > \end{gather*} 
          \end{proof} 
    
    \begin{remarks} \begin{enumerate} \item \label{Schur-remk}
    Of course the same statement holds for $S$ representing an effective ample cycle and as a consequence of Giambelli's formula, this is the case for  positive combination of Schur polynomials \footnote{And it follows from  \cite{Fulton-Lazarsfeld} that there are no other such classes. Note that
    \cite{Fulton-Lazarsfeld} is also a basic ingredient of \cite{DPS94}. }
    
    $$P=\det (c_{a_i-i+j}) \;\text{où}\;r\geq a_1\geq a_2 \geq .... \geq a_r\geq 0$$
    For example we get that for a bundle generated by its sections over a  surface,    $$ \vert \left <p_1(E_ \R)^2-p_2(E_\R), M_\R \right > \vert \leq \left < c_1(E_\C)^2-c_2(E_\C), M_\C \right > $$ 

\item Of course the same holds, with obvious changes,  if $E^*_\C$ is generated by sections, since 
then $c_j(E^*_\C)=(-1)^j c_j(E_\C)$, and $p_j(E_\R^*)=p_j(E_\R)$

\end{enumerate} 
  \end{remarks} 
   
   \begin{corollaire} [\cite{Kharlamov2}]
  Let $M^4_\R$ be real algebraic manifold, and assume $T^*M_C^4$ is  generated by sections. If $\sigma (M^4_\R )$ is the signature of the $4$-dimensional  manifold $M^4_\R$, then we have
   $$ \vert 3 \sigma (M^4_\R) \vert \leq c_2^2(M_\C)$$
   \end{corollaire} 
  \begin{proof} 
  This is obvious once we notice that $$\sigma (M_\R)=\frac{1}{3}\langle p_1(M^4_\R), M^4_\R \rangle$$
  \end{proof} 
      A similar proof yields a  result for Pontriagin classes, not only Pontriagin numbers as follows:

  \begin{proposition}
With the assumptions of the previous theorem,  let $ \mu _\R$ and $ \mu _\C$ be cohomology classes Poincar{\'e} dual to $Z_\R$ and $Z_\C$, where $Z_\C$ is a singular cycle in $M_\C$ with codimension two singularities.
  Then we have

 $$ \left\vert \left < p_{J}(E_\R)\cup \mu _\R, [M_\R]
  \right > \right\vert \leq \left <c_{J}(E_\C)^2
  \cup \mu _\C ,[M_\C]\right > $$
  \end{proposition}
 
  \begin{corollaire}\label{corpon}
  Let $E$ be a bundle of rank $r$, $ \Lambda $ be a real positive line bundle,
  and $p$ be such that $E\otimes \Lambda^p$ is globally generated.
  Then we have as a consequence of the above theorem $$c_k(E\otimes \Lambda^p)=\sum_{j=0}^{k} p^{k-j}c_1(\Lambda)^{k-j}\cdot c_j(E)$$
  Taking $p$ even, we have that $\Lambda_{\mid M_\R}$ is trivial,
  hence $p_{4r}(E_\R\otimes \Lambda_\R)=p_{4r}(E_\R)$.
  Thus  denoting $\Gamma_r(E;t)=\sum_{j=0}^{r} t^{r-j}\cdot c_j(E)$, we have

  $$ \left\vert \left <p_{J}(E_\R), [M_\R]\right > \right\vert \leq
  \left <\Gamma_{J}^2(E_\C; pc_1(\Lambda)),[M_\C]\right > $$
  \end{corollaire}

 \section{Bounds}\label{Bounds}
Let $X$ be a $n\geq 1$-dimensional connected smooth subvariety of
$\PP^m$ of degree $d$. We denote by $\L_X$ the line bundle
$\omega_X\otimes\O(n+2)$ where $\omega_X=\bigwedge^nT^*X$ is the
canonical bundle. Recall that $\L_X$ is ample. We denote by
$L_X,K_X$ and $h$ the (first) Chern classes of $\L_X,\omega_X$ and
$\mathcal {O} (1)$ respectively.

\begin{proposition}\label{hL} One has the inequalities $1\leq h^{n-i}L_X^i\leq d^{i+1}$
for $i=0,\cdots,n$.
\end{proposition}
\proof The left hand side inequality follows from the
ampleness of both $\O(1)$ and $\L_X$.

Let us prove the other inequality  by induction on $n$. For $i=0$, this
inequality reduces to the equality $h^n=d=\deg(X)$.


Assume that $n>1$. Let $H$ be a smooth (connected) hyperplane
section of $X$. It's an $n-1$-dimensional subvariety
of $\PP^{m-1}$ of  degree $d$. By adjunction, one has
$$\omega_H=\omega_X\otimes\O(1)_{|H}$$ and therefore

\begin{equation}\label{restri}
  (\L_X)_{|H}=\L_H.
\end{equation}

First, one has the equalities (projection formula)
$$h^{n-i}L_X^i=h^{n-i-1}{L_X^{i}}_{|H}=h^{n-1-i}L_H^{i}$$
for all $i<n$.

 Because both $\O(1)$ and $\L_X$ are ample, one has
$hL_X^{n-1}>0$ and therefore there exists some rational
$\alpha\in\QQ$ such that
$$(L_X+\alpha h)L_X^{n-1}=0.$$
In other words, $L_X+\alpha h$ is a primitive class in
${H^2}(X,\QQ)$. By the Hodge index theorem, one gets $$(L_X+\alpha
h)^2L_X^{n-2}\leq 0.$$ Thus  the real quadratic form
$$t\mapsto (L_X+th)^2L_X^{n-2}$$ has non negative discriminant
(because it is negative for $t=0$ and positive for large $t$). This gives the
inequality
$$0<L_X^n\leq\frac{(hL_X^{n-1})^2}{(h^2L_X^{n-2})}$$

Applying the above inequality to a smooth codimension $k$-plane
section, we get

$$0<h^kL_X^{n-k}\leq\frac{(h^{k+1}L_X^{n-k-1})2}{(h^{k+2}L_X^{n-k-2})}$$

 Set now $l_k=\log (h^{n-k}L_X^{k})$. Since $h^{n-k}L_X^{k} \geq 1$ we have $l_k\geq 0$ and
 the previous inequality implies $l_k+l_{k-2}\leq 2 l_{k-1}$, or else $l_k-l_{k-1}$ is decreasing.

 We claim that
$l_1-l_0\leq\log (d)$. Indeed, first one has certainly
$l_0=\log (h^n)=\log(d)$. Second,  observe that $l_1$ is obtained for curves as
follows. One could use Castelnuovo's inequalities, but let us just
give these elementary bounds. Let $C$ be a generic intersection of
$X$ and $(n-1)$ hyperplane. It is a genus $g$ curve of degree $d$.
By~\ref{restri}, one has
$$l_1=\log (\deg(\L_C))=\log (2g-2+3d)$$

One just has to bound $g$ in terms of $d$. A generic projection
defines a  birational morphism onto a plane degree $d$ curve $C$.
On deduces the inequality $$1-g=\chi(\O_X)\geq
\chi(\O_C)=-d(d-3)/2$$ and therefore $$l_1\leq \log (d^2)$$ giving the
inequality in this case.

 Thus $l_k\leq (k+1)\log (d)$ and this concludes our proof. \qed
\section{Chern classes} Using the splitting principle, we get the
formula
$$c_i(E\otimes L)=\sum_j\binom{n-j}{i-j}c_1(L)^{i-j}c_j(E)$$ for $E$ a
(complex) vector bundle and $L$ a line bundle on any variety.
Applying this identity to $E=\Omega(2)$ and $L=\O(-2)$

\begin{equation}\label{c(EL)}
 c_i(\Omega)=\sum_j(-2)^{i-j}\binom{n-j}{i-j}h^{i-j}c_j(\Omega(2)).
\end{equation}
Because $\Omega_{\PP^m}(2)$ is globally generated, so is its
quotient $\Omega_X(2)$ which is therefore nef. By~\cite{DPS94},
corollary 2.6, every Chern number is controlled by $K_X$ and $h$,
or, what's the same by $L_X$ and $h$. More precisely for every
multi index $I=(i_1,\cdots,i_r)$ of $r\leq n$ integers in
$[1,\cdots,n]$, one has
\begin{equation*}
  0\leq c_I(\Omega(2))h^{n-|I|}\leq
  c_1(\Omega_X(2))^{|I|}h^{n-|I|}=(L_X+(n-2)h)^{|I|}h^{n-|I|}.
\end{equation*}
where $c_I=c_{i_1}\cdots c_{i_r}$. Let's assume $n>1$ (the case
$n=1$ is left to the reader !). Using~\ref{hL}, we get the
estimate
\begin{equation}\label{cI}
0\leq c_I(\Omega(2))h^{n-|I|}\leq d(d+n-2)^{|I|}.
\end{equation}

To bound $c_I(\Omega)h^{n-|I|}$, let's denote the multi-index
$(n,\cdots,n)\in\Z^r$ by $\n$. Using~\ref{c(EL)}, one gets

\begin{align*}
  |c_I(\Omega)h^{n-|I|}| &\leq\sum_J2^{|I|-|J|}\binom{\n-J}{I-J}h^{n-|J|}c_J(\Omega(2)) \\
  & \leq\sum_J2^{|I|-|J|}\binom{\n-J}{I-J}d(d+n-2)^{|J|}\
  \text{by}\
  (\ref{cI})\\&= d(d+n-2)^{|I|}\sum_J\binom{\n-I+I-J}{\n-I}2^{|I|-|J|}(d+n-2)^{|J|-|I|}\\
\end{align*}
Applying  the identity
$$\sum_{i}\binom{i+m}{m}t^i=\frac{1}{(1-t)^{m+1}}\text{ with }
t^{-1}=(d+n-2)/2\text{ and }m=n-i_\alpha,$$ we obtain that the
last sum is bounded by
$$\prod_\alpha\frac{1}{(1-t)^{n-i_\alpha+1}}$$ which in turn is less
or equal than $2^{{n^2}}$ because $0\leq t\leq 1/2$ provided $d+n\geq 6$.  This yields
the
 estimate (the cases where $d+n\leq 6$ can be dealt with by inspection)
\begin{equation}\label{cIn}
 |c_I(\Omega)h^{n-|I|}|\leq  2^{n^2}d(d+n-2)^{|I|}
\end{equation}
where $I=(i_1,\cdots,i_r)$.

     \section{Proof of the main theorem~\ref{boundpon}}

Assume that $X$ is a real smooth subvariety of degree $d$ of some
projective space $\R P^m$. Recall that this means that $X$ is defined
by polynomials with real coefficients and that the complex
corresponding variety is smooth.  

The (twisted) cotangent bundle $\Omega_X(2)$ is therefore a
quotient of the twisted cotangent bundle $\Omega_P(2)$. Because
the latter is globally generated (straightforward computation), so
is $\Omega_X(2)$. By~\ref{corpon} and~\ref{cIn}, we get easily the
claimed inequality
  $$ \vert p_I(X(\R)) \vert \leq 2^{n^2+3n}d (d+n-2)^{n}.$$
  
  Once this is proved, the finiteness of the number of cobordism classes follows from the fact that  the cobordism ring as the work of 
  \cite{Wall}, completing the results of Thom, Rokhlin and Milnor, is determined by the Pontriagin classes.

  \begin{remark} According to remark \ref{Schur-remk}, there are also inequalities
  for the Schur's polynomials. If  $$P(c_k)=\det (c_{a_i-i+j}) \;\text{où}\;r\geq a_1\geq a_2 \geq .... \geq a_r\geq 0$$
  
  we have at least that for any Schur's polynomial $P$ there is a universal function $F_P(d,n)$ such that   $$ \vert \left < P(p_{2k}), X(\R) \right > \vert \leq F_P(d,n)$$
 \end{remark}

\section{Betti numbers} Let's explain how to obtain bounds for the
Betti numbers of $X$ using ~\ref{cIn}. It is important to note
that the results of this section only depend on the complex degree
of $M_\C$, that is lowest degree of an ample (not necessarily
real) line bundle. In fact this section is about bounding Betti
numbers of complex projective manifolds, as functions of their
degree, and the application to real manifolds is a consequence of
the Smith-Thom inequality.

If $n=1$, we have
$$h^{1,0}=g(X)\leq \frac{(d-1)(d-2)}{2}$$
which gives $$b(X)=\leq 2+(d-1)(d-2).$$

Let's assume $n>1$. Let $H$ a smooth hyperplane section as above.

By the Lefschetz Hyperplane theorem, we have the relations
$$b_i(X)= b_{i}(H)\text{ if } i<n-1\text{ and } b_{n-1}(X)\leq b_{n-1}(H).$$ By
Poincaré duality, we have
$$b_i(X)=b_{2n-i}(X).$$ Therefore, one gets
$$b_i(X)\leq b_{i-2}(H)\text{ if } i>n.$$
It remains to control the middle term, $b_n(X)$. But  the
holomorphic Gauss-Bonnet formula (due to Chern) says
$$\chi_{}(X)=c_n(T_X).$$ Notice that this formula
follows easily from the Riemann-Roch theorem and from the well
known formula
$$c_n(T_X)=\td(T_X)\sum(-1)^p\ch(\bigwedge^p(\Omega_X))$$
(\cite{Fulton-Lang}, proposition 5.3). We get therefore\footnote{Notice that the Euler characteristic is independent from the coefficient field}

\begin{align*}
  b_n(X)&\leq |\chi_{}(X)|+|\sum_{i\not=n}(-1)^ib_i(X)| \\
   &\leq |\chi_{}(X)|+|\sum_{i\not=n,n-1,n+1}(-1)^ib_i(X)|+b_{n-1}(X)+b_{n+1}(X)\\
   &=|\chi_{}(X)|+|\chi_{}(H)|+b_{n-1}(X)+b_{n+1}(X)\\
   &\leq |c_{n}(TX)|+|c_{n-1}(TH)|+2b_{n-1}(H)\\
    &\leq 2(2^{{n^2}}d^{n+1}+b_{n-1}(H))\text{ thanks to
    }(\ref{cIn})\\
\end{align*}

If $b(X)=\sum b_i(X)$ is the total Betti number, one gets
$$b(X)\leq 4b(H)+2.2^{{n^2}}d^{n+1}$$ and finally

\begin{equation*}
  b(X)\leq 4^n(d^2+1)+2\sum_{k=2}^n4^{n-k}.2^{k^2}d^{k+1}\leq
  2\sum_{k=0}^n4^{n-k}.2^{k^2}d^{k+1}
\end{equation*}

To get a bound without summation, one can for instance bound
$4^{n-k}.2^{{k^2}}$ by $2^{n+2+(k+1)(n-2)}$. The last sum is
bounded by
$$2^{n+3}\sum_{k=0}^n(2^{n-2}d)^{k+1}=2^{n+3}
\frac{(2^{n-2}d)^{n+2}-2^{n-2}d}{2^{n-2}d-1}$$ and, at last,
\begin{equation}\label{betti}
  b(X)\leq 2^{{n^2}+2}d^{n+1}
\end{equation}  (of
course it's easy to get a better bound, but with the same growth
in $d$ at fixed $n$).

Using the Smith-Thom's inequality, one gets

\begin{proposition} Let $X$ be a real subvariety of the projective
space. Assume that $X_\C$ is smooth, connected of degree $d$ and
dimension $n$. Then, one has $b(X(\R))\leq 2^{{n^2}+2}d^{n+1}$.
\end{proposition}

  Notice that this bound does not depend on the real structure. Also the bound on $\mod 2$ Betti numbers of course implies the same on rational Betti numbers. 
  
\begin{remark} It should be interesting to have some control on the
fundamental group of the (non connected) real part. For instance,
given a finite covering of $X(\R)$, can one lift this covering to
a finite ramified covering of $X_\C$ ramified over an hypersurface
of controlled degree (without real point of course) ?
\end{remark}

\section{Appendix: Some simple finiteness results}

 Let $X$ be a $n$ dimensional complex algebraic variety. We denote by $\bar X$
 the conjugate variety obtained by the base change by the conjugation $\C\ra\C$. Explicitly, if
 $X$ is (locally) defined by polynomial equations $\sum a_{i,j}x^i=0$, the variety
 $\bar X$ is defined by the conjugate equations $\sum \bar a_{i,j}x^i=0$.
 Of course, one has $\bar{\bar X}=X$. To give a real structure on $X$ remains to give a complex
 morphism $t:X\ra \bar X$ such that $\bar t\rond t=\Id_X$
 (we'll say simply that $t$ is a skew involution of $X$. The real
 points of this real complex structure we'll be denoted by
 $X_t(\R)$ (the set of fixe points).
 If $t_1,t_2$ are two skew involutions, the composite $t_1\rond\bar t_2$ belongs to $\Aut(X)$
  We get therefore the following well
 known result

 \begin{proposition} A complex variety of general type has a finite number of real structures.
 \end{proposition}
\begin{proof} Observe that $\Aut(X)$ is finite in this case
(\cite{Kobayashi}).\end{proof}

 \subsection{Final remarks} The set of morphisms $\Hom(X,\bar X)$ is parameterized by the complex points
 of a countable union of quasiprojective varieties, which is a locally closed subvariety
 of the Hilbert scheme of $X\times\bar X$ (look at the graph), and therefore the set of skew
 morphisms a the
 same property.
 Recall (Chow) that the subvariety of the Hilbert scheme of a
 smooth projective (polarized) variety parameterizing the reduced
 subvarieties of bounded degree is quasiprojective. If $L$
 is an ample line bundle on $X$, it defines a conjugate line
 bundle $\bar L$ on $\bar X$ and therefore $X\times\bar X$ is
 polarized by $L\boxtimes\bar L$. The degree of the graph of a
 morphism $t\in\Hom(X,\bar X)$ is $(c_1(L)+t^*c_1(\bar L))^n$. In
 particular, if $t^*\bar L$ and $L$ are numerically equivalent, this degree is
 bounded by $2^n\deg_L(X)$. This is certainly the case if for instance $H^{1,1}(X)=\C$.
Another good situation is when $X$ is Fano. Indeed, one can
 take $L=\omega_X^{-1}$ (observe that $\bar\omega_X=\omega_{\bar
 X}$). 

  \begin{corollaire} The real structures of a given  complex Fano variety have
  a finite number of deformation  classes.
 \end{corollaire}

\begin{remark}
If the automorphism group of $X$ is linear algebraic, its action
on the Picard group factors through the finite group of the
connected components. One deduces (averaging a given ample bundle)
that the number of deformation classes is finite in this case. A more general result using \cite{Borel-Serre} can be found in \cite{Degtyarev-Kharlamov3} (Section D.1.10). 
  
In particular, the number of diffeomorphism types of $X_t(\R)$ when
$t$ runs over all real structures is then finite. This is for
instance the case for toric varieties (due to Demazure). This last
observation was first made in \cite{Delaunay}. The same averaging
process can more generally be achieved if for instance the closure
of the ample cone is (rational) polyhedral (which is also the case
for toric varieties), giving a more elementary proof of the former
finiteness result (use the fact that the 1-dimensional edges are
permuted by the automorphisms and that a an automorphisms fixing
the edges is the identity because it is an integral dilatation on
each edge on one hand and, on the other hand, it is volume
preserving, its  determinant being equal to $\pm 1$).
\end{remark}

    \bibliography{pont-last-essai}
    \bibliographystyle{alpha} 

\end{document}